\documentclass[10pt]{article}
\usepackage{graphicx}
\usepackage{amsmath}
\usepackage[dvips]{epsfig}
\usepackage{amssymb}
\begin{document}
\begin{center} {\Large \bf  On the  generalized higher-order $q$-Bernoulli numbers and polynomials}
\\ \vspace*{12 true pt}  T.  Kim$^1$, Byungje Lee$^2$,   C. S.  Ryoo$^3$
\vspace*{12 true pt} \\
$^1$ Division of General Education-Mathematics,
 Kwangwoon University, Seoul 139-701,  Korea \\
 $^2$ Depart. of Wireless Communications Engineering,
 Kwangwoon University, Seoul 139-701,  Korea \\
              $^3$ Department of Mathematics, Hannam University,  Daejeon 306-791, Korea \\
          \end{center}
\vspace*{12 true pt} \noindent {\bf Abstract :} In this paper we
give some interesting equation of $p$-adic $q$-integrals on $
\mathbb{Z}_p $. From those $p$-adic $q$-integrals, we present a
systemic study of some families of extended Carlitz $q$-Bernoulli
numbers and polynomials in $p$-adic number field.

\vspace*{12 true pt} \noindent {\bf 2000 Mathematics Subject
Classification :} 11B68, 11S40, 11S80

 \vspace*{12 true pt} \noindent {\bf
Key words :}   Bernoulli  numbers, Bernoulli  polynomials,
generalized higher-order  $q$-Bernoulli numbers and polynomials

\begin{center} {\bf 1. Introduction } \end{center}

Let $p$ be a fixed prime number. Throughout this paper
$\mathbb{Z}_p$,  $\mathbb{Q}_p$, $\mathbb{C}$, $\mathbb{C}_p$
will, respectively,  denote the ring of $p$-adic rational integer,
the field of $p$-adic rational numbers,  the complex number field
and  the
  completion of algebraic closure of   $\mathbb{Q}_p$.
Let $ \mathbb{N} $ be the set of natural numbers and
$\mathbb{Z}_+=\{ 0 \} \cup \mathbb{N}$.

Let $\nu_p$ be the normalized exponential valuation of
$\mathbb{C}_p$ with $|p|_p=p^{-\nu_p(p)}=p^{-1}.$ When one talks
of $q$-extension, $q$ is considered as an indeterminate, a complex
number $q\in  \mathbb{C},$ or $p$-adic number $q\in\Bbb C_p .$ If
$q\in \Bbb C$,  we normally assume  $|q|<1$, and if $q\in
\mathbb{C}_p,$ we normally assume  $|1-q|_p<1$. We use the
notation
$$[x]_q =\frac{1-q^x}{1-q}.$$
The $q$-factorial is defined as
$$[n]_q !=[n]_q [n-1]_q \cdots [2]_q[1]_q$$ and  the Gaussian
$q$-binomial coefficient is defined by
$${\binom{n}{k}}_q= \dfrac{[n]_q !}{[n-k]_q![k]_q!}=\dfrac{[n]_q [n-1]_q \cdots [n-k+1]_q}{[k]_q!},  \text{ (see [9])}. \eqno(1)$$
Note that $$ \lim_{q \rightarrow
1}\binom{n}{k}_q=\binom{n}{k}=\dfrac{n(n-1)\cdots (n-k+1)}{k!}.$$
From (1), we easily see that
$$\binom{n+1}{k}_q=\binom{n}{k-1}_q+ q^k \binom{n}{k}_q= q^{n-k} \binom{n}{k-1}_q+\binom{n}{k}_q, \text{ (see [8, 11])}.$$
 For a fixed positive integer $f, (f, p)=1$, let
$$\aligned
&X= X_f= \varprojlim_N ( \mathbb{Z}/f p^N \mathbb{Z}), \quad X_1= \mathbb{Z}_p, \\
&X^*=\bigcup_{\begin{subarray}{l} 0<a<fp \\ (a,p)=1 \end{subarray}} ( a+f p \mathbb{Z}_p), \text{ and } \\
&a+fp^N \mathbb{Z}_p=\{x\in X\mid x\equiv a\pmod{fp^N}\},
\endaligned$$
where $a\in \mathbb{ Z}$  and  $0\leq a<fp^N$(see [1-14]).

We say that $f$ is a uniformly differential function at a point $a
\in \mathbb{Z}_p$ and denote this property by $f \in UD(
\mathbb{Z}_p) $ if the difference quotients
$$F_f(x,y)= \dfrac{f(x)-f(y)}{x-y}$$
have a limit $l= f'(a)$ as $(x,y) \rightarrow (a,a)$.  For $f \in
UD( \mathbb{Z}_p) $,  let us begin with the expression
$$ \dfrac{1}{[p^N]_q }\sum_{x= 0}^{ p^N-1}f(x)q^x= \sum_{0\leq x <  p^N}f(x) \mu_q(x +p^N \mathbb{Z}_p),$$
representing a $q$-analogue of the Riemann sums for  $f$, (see
[8-18]). The integral of $f$ on $ \mathbb{Z}_p$ is defined as the
limit $(N \rightarrow  \infty)$ of the sums (if exists). The
$p$-adic $q$-integral (or $q$-Volkenborn integrals of $f \in UD(
\mathbb{Z}_p) $) is defined by
$$I_q(f)=\int_{X}f(x) d\mu_q(x)= \int_{\mathbb{Z}_p}f(x)d\mu_q(x)=
\lim_{N \to \infty}\dfrac{1}{[p^N]_q }\sum_{0\leq x < p^N}f(x)q^x,
\text{ (see [12])}. \eqno(2)
$$
Carlitz's $q$-Bernoull numbers $\beta_{k, q}$ can be defined
recursively by $\beta_{0, q}=1$ and by the rule that
$$   q  \left( q \beta+1 \right)^k -   \beta_{k, q}  =\left \{\begin{array}{ll}
1, & \mbox{ if } k=1, \\
0, &  \mbox{ if } k > 1,
\end{array} \right.
 $$
with the usual convention of replacing $\beta^i$ by $\beta_{i,
q}$, (see [1-13]).

It is well known that
$$\beta_{n, q}=  \int_{\mathbb{Z}_p}[x]_q^n d\mu_q(x)=  \int_{X}[x]_q^n d\mu_q(x), n \in \mathbb{Z}_+, \text{ (see [9])},$$
and
$$\beta_{n, q}(x)=  \int_{\mathbb{Z}_p}[y+x]_q^n d\mu_q(y)=  \int_{X}[y+x]_q^n d\mu_q(y), n \in \mathbb{Z}_+,$$
where $\beta_{n, q}(x)$ are called the $n$-th Carlitz's
$q$-Bernoulli polynomials (see [9, 12, 13]).

 Let $\chi$ be the  Dirichlet's  character with
conductor $f \in \Bbb N$. Then the generalized Carlitz's
$q$-Bernoulli numbers attached to $\chi $ are defined as follows:
$$\beta_{n,\chi,  q}=  \int_{X} \chi(x) [x]_q^n d\mu_q(x), \text{ (see [13])}.$$
Recently, many authors have studied in the different several areas
related to $q$-theory (see [1-13]).  In this paper we present a
systemic study of some families of multiple Carlitz's
$q$-Bernoulli numbers and polynomials by using the integral
equations of $p$-adic $q$-integrals on $\mathbb{Z}_p$. First, we
derive some interesting the equations of $p$-adic $q$-integrals on
$\mathbb{Z}_p$. From these equations, we give some interesting
formulae for the higher-order Carlitz's $q$-Bernoulli numbers and
polynomials in the $p$-adic number field.

\bigskip
\begin{center} {\bf 2. On the generalized  higher-order $q$-Bernoulli numbers and polynomials} \end{center}
\bigskip

In this section  we assume that  $q\in \mathbb{C}_p$ with
$|1-q|_p<1$. We first consider the $q$-extension of Bernoulli
polynomials as follows:

$$ \sum_{n=0}^\infty \beta_{n, q}(x) \dfrac{t^n}{n!}= \int_{\mathbb{Z}_p} q^{-y} e^{[x+y]_q  t}
d\mu_q(t)=-t \sum_{m=0}^\infty  e^{[x+m]_q  t} q^{x+m}. \eqno(3)$$
From (3), we note that
$$ \aligned   \beta_{n, q}(x) &= \dfrac{1}{(1-q)^n} \sum_{l=0}^n \binom nl (-q^x)^l \dfrac{l}{[l]_q} \\
&=\dfrac{1}{(1-q)^{n-1}} \sum_{l=0}^n \binom nl (-q^x)^l \left( \dfrac{l}{1-q^l} \right)\\
&=\dfrac{n}{(1-q)^{n-1}} \sum_{l=0}^{n-1} \binom {n-1}{l}
q^{(l+1)x} \left(  \dfrac{1}{1-q^{l+1}} \right) (-1)^{l+1} \\
&=   \dfrac{-n}{(1-q)^{n-1}} \sum_{m=0}^{\infty}  q^{m+x}
\sum_{l=0}^{n-1}  \binom {n-1}{l} q^{l(x+m)} \\
&= -n  \sum_{m=0}^{\infty}  q^{m+x} [x+m]_q^{n-1}.\endaligned
\eqno(4)
 $$
Note that $$ \lim_{ q \rightarrow 1} \beta_{n,q}(x)=-n
\sum_{m=0}^{\infty}  (x+m)^{n-1}=B_n(x),$$ where $B_n(x)$ are
called the $n$-th ordinary Bernoulli polynomials. In the special
case,  $x=0$, $\beta_{n, q}(0)= \beta_{n,  q }$ are called the
$n$-th  $q$-Bernoulli numbers.

By (4), we have the following lemma.

\medskip
{ \bf Lemma 1.} For $ n  \geq 0,$ we have
$$ \aligned     \beta_{n, q}(x)  & =  \int_{\mathbb{Z}_p} q^{-y} [x+y]_q^n  d\mu_q(y) = -n  \sum_{m=0}^{\infty}  q^{m+x} [x+m]_q^{n-1} \\
&= \dfrac{1}{(1-q)^n} \sum_{l=0}^n \binom nl (-q^x)^l
\dfrac{l}{[l]_q}. \endaligned
 $$
\medskip

Now, we consider  the $q$-Bernoulli polynomials of order  $r \in
\mathbb{N}$ as follows:

$$ \sum_{n=0}^{ \infty}  \beta_{n, q}^{(r)}(x) \dfrac{t^n}{n!}=
 \underset{r \mbox{ times}}{\underbrace{ \int_{ \mathbb{Z}_p }\cdots \int_{ \mathbb{Z}_p} } }
     q^{-(x_1+ \cdots + x_r)} e^{[x+ x_1+ \cdots + x_r ]_q t}    d\mu_{q}(x_1)  \cdots d\mu_{q}(x_r). \eqno(5)$$
By (5), we see that
$$ \aligned  \beta_{n, q}^{(r)}(x) & =
 \underset{r \mbox{ times}}{\underbrace{ \int_{ \mathbb{Z}_p }\cdots \int_{ \mathbb{Z}_p} } }
     q^{-(x_1+ \cdots + x_r)} [x+ x_1+ \cdots + x_r ]_q^n    d\mu_{q}(x_1)  \cdots d\mu_{q}(x_r)\\
     &= \dfrac{1}{(1-q)^n} \sum_{l=0}^n \binom nl (-1)^l q^{xl}
\left( \dfrac{l}{[l]_q} \right)^r. \endaligned $$
 In the special
case,  $x=0$, the sequence  $\beta_{n, q}^{(r)}(0)= \beta_{n, q
}^{(r)}$ is refereed as the  $q$-extension of Bernoulli numbers of
order $r$. For $f \in \mathbb{N}$, we have
$$ \aligned  \beta_{n, q}^{(r)}(x) & =
 \underset{r \mbox{ times}}{\underbrace{ \int_{ X }\cdots \int_{ X }
 }}     q^{-(x_1+ \cdots + x_r)} [x+ x_1+ \cdots + x_r ]_q^n    d\mu_{q}(x_1)  \cdots d\mu_{q}(x_r)\\
     &= \dfrac{1}{(1-q)^n} \sum_{l=0}^n \dfrac{ \binom nl (-1)^l
     q^{l(x+ a_1 +\cdots + a_r )}l^r }{[lf]_q^r} \\
     &= [f]_q^{n-r}  \sum_{a_1, \cdots, a_r =0}^{f-1}   \beta_{n, q^f}^{(r)}\left( \dfrac{a_1+\cdots + a_r+x }{f}\right). \endaligned  \eqno(6)$$
By (5) and (6), we obtain the following theorem.

\medskip
{ \bf Theorem 2.} For $r \in \mathbb{Z}_+, f \in \Bbb N, $ we have
$$ \aligned  \beta_{n, q}^{(r)}(x)
& = \dfrac{1}{(1-q)^n} \sum_{l=0}^n   \sum_{a_1, \cdots, a_r =0}^{f-1}  \binom nl (-1)^l q^{ l ( a_1+\cdots +a_r +x)} \dfrac{l^r }{[lf]_q^r} \\
     &= [f]_q^{n-r}  \sum_{a_1, \cdots, a_r =0}^{f-1}   \beta_{n, q^f}^{(r)}\left( \dfrac{a_1+\cdots + a_r+x }{f}\right). \endaligned $$
\medskip

Let $\chi$ be the  primitive Dirichlet's  character with conductor
$f \in \Bbb N$. Then the generalized $q$-Bernoulli polynomials
attached to $\chi$ are defined by

$$ \sum_{n=0}^{ \infty}  \beta_{n, \chi, q}(x) \dfrac{t^n}{n!}=
  \int_{ X }
     \chi(y)  q^{-y} e^{[x+ y ]_q t}   d\mu_{q}(y) . \eqno(7) $$
From  (7),  we derive
$$ \aligned   \beta_{n, \chi, q}(x) &= \int_{ X }
     \chi(y)  q^{-y} {[x+ y ]_q ^n}   d\mu_{q}(y) \\
     &= \sum_{a=0}^{f-1} \chi(a) \lim_{N \rightarrow \infty} \dfrac{1}{[f p^N ]_q} \sum_{y=0}^{fp^N-1}[a+x+fy]_q^n\\
&=\dfrac{1}{(1-q)^n} \sum_{a=0}^{f-1} \chi(a)  \sum_{l=0}^n   \binom nl (-1)^l q^{ l ( x+a)} \dfrac{l }{[lf]_q} \\
&=\sum_{a=0}^{f-1} \chi(a)  \sum_{m=0}^{\infty} \left( -n [x+a+mf]_q^{n-1} \right)  \\
&= -n  \sum_{m=0}^{\infty} \chi(m) [x+m]_q^{n-1}.\endaligned
\eqno(8)
 $$
By (7) and (8), we can give the generating function for the
generalized $q$-Bernoulli polynomials attached to $\chi$ as
follows:

$$ F_{ \chi, q }( x, t )= - t  \sum_{m =0}^{
\infty}   \chi( m )  e^{[x+m]_q t}=\sum_{n=0}^\infty \beta_{n,
\chi, q}(x) \dfrac{t^n}{n!} . \eqno(9)$$

From (1), (8) and (9), we note that

$$ \aligned   \beta_{n, \chi, q}(x) & = \dfrac{1}{ [f]_q} \sum_{a=0}^{f-1} \chi(a)   \int_{ \mathbb{Z}_p }
      q^{-fy} [a+x+fy]_q^n   d\mu_{q^f}(y) \\
     &= [f]_q^{n-1} \sum_{a=0}^{f-1} \chi(a) \beta_{n, q^f} \left( \dfrac{a+x}{f} \right).\endaligned
 $$
In the special case,   $x=0$, the sequence  $\beta_{n, \chi,
q}(0)= \beta_{n, \chi, q}$ are called the $n$-th generalized
$q$-Bernoulli numbers
 attaches to $\chi$.

Let us consider the  higher-order $q$-Bernoulli polynomials
attached to $\chi$ as follows:
$$  \underset{r \mbox{ times}}{\underbrace{ \int_{X }\cdots \int_{ X} } }
      \left(  \prod_{i=1}^r \chi( x_i) \right) e^{[x+ x_1+ \cdots + x_r ]_q t}  q^{-(x_1+ \cdots + x_r)}   d\mu_{q}(x_1)  \cdots d\mu_{q}(x_r)
     =\sum_{n=0}^{ \infty}  \beta_{n, \chi, q}^{(r)}(x) \dfrac{t^n}{n!}, \eqno(10)$$
where $\beta_{n, \chi, q}^{(r)}(x) $ are called the $n$-th
generalized $q$-Bernoulli polynomials of order $r$
 attaches to $\chi$.

By (10), we see that

$$ \aligned   \beta_{n, \chi, q}^{(r)}(x) &=
\dfrac{1}{(1-q)^n}  \sum_{l=0}^n   \binom nl (-q^x)^l  \sum_{a_1,
\cdots, a_r =0}^{ f-1}  \left(  \prod_{i=1}^r \chi( a_i) \right)  q^{ l  \sum_{i=1}^r a_i } \dfrac{l^r }{[lf]_q^r} \\
&= [f]_q^{n-r}  \sum_{a_1, \cdots, a_r =0}^{f-1}  \left(
\prod_{i=1}^r \chi( a_i) \right)  \beta_{n, q^f}^{(r)} \left(
\dfrac{x+a_1+\cdots + a_r }{f}\right).
\endaligned  \eqno(11)
 $$
In the special case,   $x=0$,  the sequence  $\beta_{n, \chi,
q}^{(r)}(0)= \beta_{n, \chi, q}^{(r)}$ are called the $n$-th
generalized $q$-Bernoulli numbers of order $r$
 attaches to $\chi$.

By (10) and (11), we obtain the following theorem.

\medskip
{ \bf Theorem 3.} Let $\chi$ be the  primitive Dirichlet's
character with conductor $f \in \Bbb N$.  For $ n \in \Bbb Z_+, r
\in \mathbb{N},$  we have
$$ \aligned   \beta_{n, \chi, q}^{(r)}(x) &=
\dfrac{1}{(1-q)^n}  \sum_{l=0}^n   \binom nl (-q^x)^l  \sum_{a_1,
\cdots, a_r =0}^{ f-1}  \left(  \prod_{i=1}^r \chi( a_i) \right)  q^{ l  \sum_{i=1}^r a_i } \dfrac{l^r }{[lf]_q^r} \\
&= [f]_q^{n-r}  \sum_{a_1, \cdots, a_r =0}^{f-1}  \left(
\prod_{i=1}^r \chi( a_i) \right)  \beta_{n, q^f}^{(r)} \left(
\dfrac{x+a_1+\cdots + a_r }{f}\right).
\endaligned
 $$
\medskip

 For $ h \in \Bbb Z$, and $ r
\in \mathbb{N},$  we introduce the extended higher-order
$q$-Bernoulli polynomials as follows:
$$  \beta_{n, q}^{(h, r)}(x)  =
 \underset{r \mbox{ times}}{\underbrace{ \int_{ \mathbb{Z}_p }\cdots \int_{ \mathbb{Z}_p} } }
     q^{ \sum_{j=1}^r (h-j-1)x_j} [x+ x_1+ \cdots + x_r ]_q^n    d\mu_{q}(x_1)  \cdots
     d\mu_{q}(x_r). \eqno(12)$$
From (12), we note that
$$  \beta_{n, q}^{(h, r)}(x)  =
  \dfrac{1}{(1-q)^n} \sum_{l=0}^n \binom nl (-1)^l q^{lx}
  \dfrac{\binom{l+h-1}{r}}{ \binom{l+h-1}{r}_q }
  \dfrac{r!}{[r]_q!},
 \eqno(13)$$
 and
$$\beta_{n, q}^{(h, r)}(x)  = [f]_q^{n-r}  \sum_{a_1, \cdots, a_r =0}^{f-1} q^{ \sum_{j=1}^r (h-j)a_j}
\beta_{n, q^f}^{(h,r)} \left(
\dfrac{x+a_1+\cdots + a_r }{f}\right).$$

 In the special
case,   $x=0$, $\beta_{n, q}^{(h, r)}(0)= \beta_{n, q }^{(h, r)}$
are called  the  $n$-th $(h, q)$-Bernoulli numbers of order $r$.

By (13), we obtain the following theorem.

\bigskip
 { \bf Theorem 4.}  For  $ h \in \Bbb Z, r \in \mathbb{N},$
we have
$$  \beta_{n, q}^{(h, r)}(x)  =
  \dfrac{1}{(1-q)^n} \sum_{l=0}^n \binom nl (-q^x)^{l}
  \dfrac{\binom{l+h-1}{r}}{ \binom{l+h-1}{r}_q }
  \dfrac{r!}{[r]_q!},
 $$
 and
$$\beta_{n, q}^{(h, r)}(x)  = [f]_q^{n-r}  \sum_{a_1, \cdots, a_r =0}^{f-1} q^{ \sum_{j=1}^r (h-j)a_j}
\beta_{n, q^f}^{(h,r)} \left( \dfrac{x+a_1+\cdots + a_r
}{f}\right).$$

\bigskip

Let $\chi$ be the  primitive Dirichlet's  character with conductor
$f \in \Bbb N$. Then we consider  the generalized $(h,
q)$-Bernoulli polynomials attached to $\chi$  of order $r$ as
follows:
$$ \beta_{n, \chi, q}^{(h, r)}(x)  =
 \underset{r \mbox{ times}}{\underbrace{ \int_{ X }\cdots \int_{ X }
 }}    q^{ \sum_{j=1}^r (h-j-1)x_j} \left(
\prod_{j=1}^r \chi( x_j) \right) [x+ x_1+ \cdots + x_r ]_q^n
d\mu_{q}(x_1)  \cdots d\mu_{q}(x_r). \eqno(14)
$$
By (14), we see that

$$ \beta_{n, \chi, q}^{(h, r)}(x)  =
 [f]_q^{n-r}  \sum_{a_1, \cdots, a_r =0}^{f-1} q^{ \sum_{j=1}^r (h-j)a_j} \left(
\prod_{j=1}^r \chi( a_j) \right) \beta_{n, q^f}^{(h,r)}\left(
\dfrac{x+ a_1+\cdots + a_r }{f}\right).
  \eqno(15)$$
In the special case,   $x=0$, $\beta_{n, \chi, q}^{(h, r)}(0)=
\beta_{n, \chi,  q }^{(h, r)}$ are called  the  $n$-th generalized
$(h, q)$-Bernoulli numbers  attached to $\chi$ of order $r$.

From (14) and (15), we note that
$$\beta_{n, \chi,  q }^{(h, r)}=(q-1)\beta_{n+1, \chi,  q }^{(h-1, r)}+\beta_{n, \chi,  q }^{(h-1, r)}.$$
By (12), it is easy to show that
$$ \aligned   & \beta_{n, \chi, q}^{(h, r)} =
\underset{r \mbox{ times}}{\underbrace{ \int_{ \mathbb{Z}_p
}\cdots \int_{ \mathbb{Z}_p }
 }}    [ x_1+ \cdots + x_r ]_q^n q^{ \sum_{j=1}^r (h-j-1)x_j}
d\mu_{q}(x_1)  \cdots d\mu_{q}(x_r) \\
&=  \int_{ \mathbb{Z}_p }\cdots \int_{ \mathbb{Z}_p }
    [ x_1+ \cdots + x_r ]_q^n  \{ [x_1+ \cdots + x_r]_q(q-1)+1 \} q^{ \sum_{j=1}^r (h-j-2)x_j}
d\mu_{q}(x_1)  \cdots d\mu_{q}(x_r).
\endaligned
\eqno(16)
 $$
Thus, we have
$$\beta_{n, q}^{(h, r)} = (q-1)\beta_{n+1, q}^{(h-1, r)} + \beta_{n,  q}^{(h-1, r)}. $$
From (12) and (16), we can also derive

$$ \aligned   &  \underset{r \mbox{ times}}{\underbrace{ \int_{ \mathbb{Z}_p
}\cdots \int_{ \mathbb{Z}_p }
 }}    q^{(n-2)x_1+ (n-3)x_2+ \cdots+ (n-r-1)x_r} d\mu_{q}(x_1)  \cdots d\mu_{q}(x_r) \\
&=  \int_{ \mathbb{Z}_p }\cdots \int_{ \mathbb{Z}_p }
  q^{-(x_1+  \cdots+ x_r)} q^{n(x_1+ \cdots+ x_r)} q^{-x_1-2x_2 - \cdots - r x_r} d\mu_{q}(x_1)  \cdots
 d\mu_{q}(x_r)\\
 &= \sum_{l=0}^n \binom nl (q-1)^l \int_{ \mathbb{Z}_p
}\cdots \int_{ \mathbb{Z}_p }
 [x_1+ \cdots + x_r]_q^l q^{-(x_1+  \cdots+ x_r)} q^{-x_1-2x_2 - \cdots - r x_r} d\mu_{q}(x_1)  \cdots
 d\mu_{q}(x_r) \\
 &= \sum_{l=0}^n \binom nl (q-1)^l \beta_{l, q}^{(0, r)},
\endaligned
\eqno(17)
 $$
and
$$  \int_{ \mathbb{Z}_p}\cdots \int_{ \mathbb{Z}_p }
   q^{(n-2)x_1+ (n-3)x_2+ \cdots+ (n-r-1)x_r} d\mu_{q}(x_1)  \cdots d\mu_{q}(x_r) = \dfrac{\binom{n-1}{r}}{\binom{n-1}{r}_q}
 \dfrac{r!}{[r]_q!}. \eqno(18)$$

It is easy to see that
$$   \sum_{j=0}^n \binom nj (q-1)^j \int_{ \mathbb{Z}_p}
 [x]_q^j q^{(h-2)x} d\mu_{q}(x)  = \int_{ \mathbb{Z}_p} ((q-1) [x]_q +1)^n  q^{(h-2)x}
 d\mu_{q}(x)= \dfrac{n+h-1}{[n+h-1]_q}.
\eqno(19)
 $$

By (16), (17), (18) and (19), we obtain the following theorem.

\medskip
{ \bf Theorem 5.} For $ h \in \Bbb Z, r \in \mathbb{N}$  and $n
\in \Bbb Z_+,$ we have
$$ \beta_{n, q}^{(h, r)} = (q-1)\beta_{n+1, q}^{(h-1, r)} + \beta_{n,  q}^{(h-1,
r)},$$
 and
$$   \sum_{l=0}^n \binom nl (q-1)^l \beta_{l, q}^{(0, r)}=
\dfrac{\binom{n-1}{r}}{\binom{n-1}{r}_q}
 \dfrac{r!}{[r]_q!}.$$
 Furthermore, we get
 $$ \sum_{l=0}^n \binom nl (q-1)^l \beta_{l, q}^{(h, 1)}=\dfrac{n+h-1}{[n+h-1]_q}. $$
\medskip

Now, we consider the polynomials of $ \beta_{n, q}^{(0, r)}(x)$ by
$$ \aligned  \beta_{n, q}^{(0, r)}(x) & =
 \underset{r \mbox{ times}}{\underbrace{ \int_{  \mathbb{Z}_p }\cdots \int_{  \mathbb{Z}_p }
 }}    [x+ x_1+ \cdots + x_r ]_q^n   q^{-2x_1-3x_2- \cdots -(r-1)x_r}  d\mu_{q}(x_1)  \cdots d\mu_{q}(x_r)\\
     &= \dfrac{1}{(1-q)^n} \sum_{l=0}^n  \binom nl (-1)^l
     q^{lx} \dfrac{\binom{l-1}{r}}{\binom{l-1}{r}_q}
 \dfrac{r!}{[r]_q!}. \endaligned  \eqno(20)
 $$
By (20), we obtain the following theorem.

\medskip
{ \bf Theorem 6.} For $ r \in \mathbb{N}$  and $n \in \Bbb Z_+,$
we have
$$ (1-q)^n \beta_{n, q}^{(0, r)}(x)  = \sum_{l=0}^n  \binom nl (-1)^l
     q^{lx} \dfrac{\binom{l-1}{r}}{\binom{l-1}{r}_q}
 \dfrac{r!}{[r]_q!}.$$
\medskip

By using multivariate $p$-adic $q$-integral on $ \Bbb Z_p$, we see
that
$$ \aligned  & q^{nx} \dfrac{\binom{n-1}{r}}{\binom{n-1}{r}_q}
 \dfrac{r!}{[r]_q!}  \\
& =  \underset{r \mbox{ times}}{\underbrace{ \int_{ \mathbb{Z}_p
}\cdots \int_{ \mathbb{Z}_p } }}  q^{nx+ (n-2)x_1+ \cdots +
(n-r-1)x_r} d\mu_{q}(x_1)  \cdots d\mu_{q}(x_r) \\
 & =  \int_{ \mathbb{Z}_p }\cdots \int_{ \mathbb{Z}_p }   \left( (q-1)[x+ x_1+\cdots +
x_r]_q+1 \right)^n
 q^{-2x_1 \cdots -(r+1)x_r} d\mu_{q}(x_1)  \cdots d\mu_{q}(x_r) \\
 &= \sum_{l=0}^n \binom nl (q-1)^l \int_{ \mathbb{Z}_p
}\cdots \int_{ \mathbb{Z}_p } [x+x_1+ \cdots + x_r]_q^l  q^{-2x_1 \cdots -(r+1)x_r} d\mu_{q}(x_1)  \cdots d\mu_{q}(x_r) \\
 &= \sum_{l=0}^n \binom nl (q-1)^l \beta_{l, q}^{(0, r)}(x).
\endaligned
 $$
Therefore, we obtain the following corollary.

\medskip
{ \bf Corollary 7.} For $ r \in \mathbb{N}$  and $n \in \Bbb Z_+,$
we have
$$ q^{nx} \dfrac{\binom{n-1}{r}}{\binom{n-1}{r}_q}
 \dfrac{r!}{[r]_q!}   = \sum_{l=0}^n  \binom nl (q-1)^l
     \beta_{l, q}^{(0, r)}(x).$$
\medskip

It is easy to show that
$$ \aligned  & \underset{r \mbox{ times}}{\underbrace{ \int_{ \mathbb{Z}_p
}\cdots \int_{ \mathbb{Z}_p } }}  [x+x_1+ \cdots + x_r]_q^n
q^{-2x_1 \cdots -(r+1)x_r} d\mu_{q}(x_1)  \cdots d\mu_{q}(x_r) \\
&=[f]_q^{n-r} \sum_{ i_1, \cdots, i_r=0}^{f-1} q^{- \sum_{l=1}^r l
 i_l} \\
 & \qquad \qquad \times     \int_{ \mathbb{Z}_p }\cdots \int_{ \mathbb{Z}_p }
 q^{- f \sum_{l=1}^r (l+1) x_l}
 \left[ \dfrac{x+ \sum_{l=1}^r i_l }{f}+ \sum_{l=1}^r x_l \right]_{q^f}^n d\mu_{q^f}(x_1)  \cdots d\mu_{q^f}(x_r) .
\endaligned
\eqno(21)
 $$
From (21), we note that

$$ \beta_{n, q}^{(0, r)}(x)  =
 [f]_q^{n-r}  \sum_{i_1, \cdots, i_r =0}^{f-1} q^{ -i_1-2i_2- \cdots -  r i_r} \beta_{n, q^f}^{(0,r)}\left(
\dfrac{x+ i_1+\cdots + i_r }{f}\right).
  $$
From the  multivariate $p$-adic $q$-integral on $ \Bbb Z_p$, we
have
$$ \aligned  & \underset{r \mbox{ times}}{\underbrace{ \int_{ \mathbb{Z}_p
}\cdots \int_{ \mathbb{Z}_p } }}  \left[ x+ x_1+\cdots +
x_r \right]_q^n  q^{-2x_1 -3x_2- \cdots -(r+1)x_r} d\mu_{q}(x_1)  \cdots d\mu_{q}(x_r)  \\
& =   \int_{ \mathbb{Z}_p }\cdots \int_{ \mathbb{Z}_p }   \left(
[x]_q + q^x[ x_1+\cdots + x_r]_q \right)^n q^{-2x_1 -3x_2- \cdots -(r+1)x_r}  d\mu_{q}(x_1)  \cdots d\mu_{q}(x_r) \\
 &= \sum_{l=0}^n \binom nl [x]_q^{n-l} q^{lx}  \int_{ \mathbb{Z}_p
}\cdots \int_{ \mathbb{Z}_p } [ x_1+ \cdots + x_r]_q^l  q^{-2x_1
-3x_2- \cdots -(r+1)x_r}  d\mu_{q}(x_1)  \cdots d\mu_{q}(x_r) ,
\endaligned
\eqno(22)
 $$
and
$$ \aligned  & \underset{r \mbox{ times}}{\underbrace{ \int_{ \mathbb{Z}_p
}\cdots \int_{ \mathbb{Z}_p } }}  \left[ x+ y+ x_1+\cdots +
x_r \right]_q^n  q^{-2x_1 -3x_2- \cdots -(r+1)x_r} d\mu_{q}(x_1)  \cdots d\mu_{q}(x_r)  \\
 &= \sum_{l=0}^n \binom nl [y]_q^{n-l} q^{ly}  \int_{ \mathbb{Z}_p
}\cdots \int_{ \mathbb{Z}_p } [x+ x_1+ \cdots + x_r]_q^l  q^{-2x_1
-3x_2- \cdots -(r+1)x_r}  d\mu_{q}(x_1)  \cdots d\mu_{q}(x_r) .
\endaligned
\eqno(23)
$$
By (22) and (23), we obtain the following corollary.

\medskip
{ \bf Corollary 8.} For $ r \in \mathbb{N}$  and $n \in \Bbb Z_+,$
we have
$$ \beta_{n, q}^{(0, r)}(x) =  \sum_{l=0}^n  \binom nl [x]_q^{n-l}
q^{lx}      \beta_{l, q}^{(0, r)},$$ and
$$\beta_{n, q}^{(0, r)}(x+y) =  \sum_{l=0}^n  \binom nl [y]_q^{n-l}
q^{ly}      \beta_{l, q}^{(0, r)}(x).$$

\medskip
Now, we also consider the polynomial of $ \beta_{n, q}^{(h,
1)}(x)$. From the integral equation on $ \Bbb Z_p$, we note that
$$ \aligned  \beta_{n, q}^{(h, 1)}(x)  & =  \int_{ \mathbb{Z}_p
}  \left[ x+ x_1 \right]_q^n  q^{x_1(h-2)}  d\mu_{q}(x_1)   \\
 &= \dfrac{1}{(1-q)^n} \sum_{l=0}^n  \binom nl (-1)^l
     q^{lx} \dfrac{l+h-1}{[l+h-1]_q}. \endaligned  \eqno(24)
 $$
By (24), we easily get
$$ \aligned  \beta_{n, q}^{(h, 1)}(x)  & =  \dfrac{1}{(1-q)^{n-1}} \sum_{l=0}^n \dfrac{\binom nl (-1)^l
     q^{lx} l }{1-q^{l+h-1}} + \dfrac{h-1}{(1-q)^{n-1}} \sum_{l=0}^n \dfrac{\binom nl (-1)^l
     q^{lx}  }{1-q^{l+h-1}} \\
     &=  \dfrac{-n}{(1-q)^{n-1}} \sum_{l=0}^{n-1} \dfrac{\binom {n-1}{l}(-1)^l q^x q^{lx} }{1-q^{l+h}}
     + \dfrac{h-1}{(1-q)^{n-1}} \sum_{l=0}^n \dfrac{\binom nl (-1)^l
     q^{lx}  }{1-q^{l+h-1}}\\
     &= -n \sum_{m=0}^\infty q^{hm+x} [x+m]_q^{n-1} +(h-1)(1-q) \sum_{m=0}^\infty q^{(h-1)m} [x+m]_q^{n}. \endaligned
 $$
Thus, we obtain the following theorem.

\medskip
{ \bf Theorem 9.} For $ h \in \mathbb{Z}$  and $n \in \Bbb Z_+,$
we have
$$ \beta_{n, q}^{(h, 1)}(x)
=  -n \sum_{m=0}^\infty q^{hm+x} [x+m]_q^{n-1} +(h-1)(1-q) \sum_{m=0}^\infty q^{(h-1)m} [x+m]_q^{n}.$$
\medskip

From the definition of  $p$-adic $q$-integral on $ \Bbb Z_p$, we
note that
$$  \int_{ \mathbb{Z}_p
}   q^{(h-2) x_1} [ x+ x_1 ]_q^n   d\mu_{q}(x_1)   =
\dfrac{1}{[f]_q} \sum_{i=0}^{f-1}  q^{(h-1)i} [i]_q^n \int_{
\mathbb{Z}_p} \left[ \dfrac{x+i}{f}+x_1 \right]_{q^f}^n
q^{f(h-2)x_1} d\mu_{q^f}(x_1) .
 $$
 Thus, we have

 $$\beta_{n, q}^{(h,1)}(x)= \dfrac{1}{[f]_q}
 \sum_{i=0}^{f-1}  q^{(h-1)i} [i]_q^n \beta_{n, q^f}^{(h,1)}\left( \dfrac{x+i}{f} \right).$$

By (24), we easily get
$$   \int_{ \mathbb{Z}_p
}  \left[ x+ x_1 \right]_q^n  q^{x_1(h-2)}  d\mu_{q}(x_1)= q^{-x}
\int_{ \mathbb{Z}_p } [x+x_1]_q^n \{  [x+x_1]_q (q-1)+1 \}
q^{x_1(h-3)}d\mu_{q}(x_1). \eqno(25)
$$
From (25), we have
$$\beta_{n, q }^{(h, 1)}(x) =q^{-x} \left( (q-1) \beta_{n+1,q }^{(h-1, 1)}(x)+ \beta_{n,  q }^{(h-1, 1)}(x) \right).$$
That is,
$$ q^x \beta_{n, q }^{(h, 1)}(x) = (q-1) \beta_{n+1,q }^{(h-1, 1)}(x)+ \beta_{n,  q }^{(h-1, 1)}(x).$$

By (24) and (25), we easily see that
$$ \aligned \int_{ \mathbb{Z}_p
}   q^{(h-2) x_1} [ x+ x_1 ]_q^n   d\mu_{q}(x_1)  & =  \int_{
\mathbb{Z}_p}   q^{(h-2) x_1} ( [ x]_q + q^x[ x_1 ]_q)^n   d\mu_{q}(x_1)  \\
     & = \sum_{l=0}^n  \binom nl [x]_q^{n-l} q^{lx} \int_{ \mathbb{Z}_p}
       q^{(h-2) x_1} [ x_1 ]_q^l   d\mu_{q}(x_1), \endaligned
 \eqno(26)
 $$
 and
$$ \aligned   & q^{h-1} \int_{ \mathbb{Z}_p
}   q^{(h-2) x_1} [ x_1+1+x ]_q^n   d\mu_{q}(x_1)  -   \int_{
\mathbb{Z}_p}   q^{(h-2) x_1}  [ x+  x_1 ]_q^n   d\mu_{q}(x_1)  \\
     & =q^x n [x]_q^{n-1} + h (q-1) [x]_q^{n} -(q-1) [x]_q^{n}. \endaligned
 $$
 For $x=0$, this give
$$   q^{h-1} \int_{ \mathbb{Z}_p
}   q^{(h-2) x_1} [ x_1+1 ]_q^n   d\mu_{q}(x_1)  -   \int_{
\mathbb{Z}_p}   q^{(h-2) x_1}  [  x_1 ]_q^n   d\mu_{q}(x_1) =\left
\{\begin{array}{ll}
1, & \mbox{ if } n=1, \\
0, &  \mbox{ if } n > 1,
\end{array} \right.
 \eqno(27) $$
 and
 $$ \beta_{0, q }^{(h, 1)}= \int_{ \mathbb{Z}_p
}   q^{(h-2) x_1}   d\mu_{q}(x_1) =\dfrac{h-1}{[h-1]_q}.
$$
From (26) and (27), we can derive the recurrence relation for
$\beta_{n, q }^{(h, 1)}$ as follows:
$$q^{h-1}\beta_{n, q }^{(h, 1)}(1)-\beta_{n, q }^{(h, 1)}=\delta_{n, 1}, \eqno(28) $$
where $\delta_{n,1}$ is kronecker symbol.

By (26), (27) and (28), we obtain the following theorem.

\bigskip
{ \bf Theorem 10.} For $ h \in \mathbb{Z}$   and $n \in \Bbb Z_+,$
we have
$$ \beta_{n, q}^{(h, 1)}(x)
=   \sum_{l=0}^n \binom nl [x]_q^{n-l} q^{lx}\beta_{l, q}^{(h,
1)}, $$ and
$$ q^{h-1}\beta_{n, q }^{(h, 1)}(x+1)-\beta_{n, q }^{(h, 1)}=q^x n [x]_q^{n-1} + h(q-1)[x]_q^n -(q-1)[x]_q^n.$$
Furthermore,
$$q^{h-2}(q-1) \beta_{n+1, q }^{(h-1, 1)}(1)+ q^{h-2}\beta_{n, q }^{(h-1, 1)}(1)-\beta_{n, q }^{(h, 1)}=\delta_{n, 1}, $$
where $\delta_{n,1}$ is kronecker symbol.

\bigskip

From the definition of  $p$-adic $q$-integral on $ \Bbb Z_p$, we
note that
$$  \int_{ \mathbb{Z}_p}   q^{-(h-2) x_1} [1-x+ x_1 ]_{q^{-1}}^n   d\mu_{q^{-1}}(x_1)
= (-1)^n q^{n+h-2} \int_{ \mathbb{Z}_p }   q^{(h-2) x_1} [x+ x_1
]_{q}^n d\mu_{q}(x_1). \eqno(29)
 $$
 By (29), we see that
$$ \beta_{n, q^{-1}}^{(h, 1)}(1-x)= (-1)^n q^{n+h-2} \beta_{n, q}^{(h,
1)}(x).$$
Note that
$$B_n(1-x)= \lim_{q \rightarrow 1} \beta_{n, q^{-1}}^{(h, 1)}(1-x)=\lim_{q \rightarrow 1}(-1)^n q^{n+h-2} \beta_{n, q}^{(h,
1)}(x)=(-1)^n B_n(x), $$ where $B_n(x)$ are the $n$-th ordinary
Bernoulli polynomials.

In the special case, $x=1$, we get
$$ \beta_{n, q^{-1}}^{(h, 1)}=(-1)^n q^{n+h-2} \beta_{n, q}^{(h,
1)}(1)= (-1)^n q^{n-1} \beta_{n, q}^{(h, 1)} \text { if } n >1. $$
It is not difficult to show that
$$ [f]_q^{n-1}  \sum_{l=0}^{f-1} q^{l(h-1)}
  \int_{ \mathbb{Z}_p } \left[ x+ \dfrac{l}{f}+ x_1 \right]_{q^f}^n
    q^{f(h-2) x_1}  d\mu_{q^{f}}(x_1)= \int_{ \mathbb{Z}_p } \left[ fx+  x_1 \right]_{q}^n
    q^{(h-2) x_1}  d\mu_{q}(x_1), f \in \Bbb N.  $$
That is,
$$ [f]_q^{n-1}  \sum_{l=0}^{f-1} q^{l(h-1)}
  \beta_{n, q^{f}}^{(h, 1)} \left( x+ \dfrac{l}{f} \right) = \beta_{n, q}^{(h, 1)}(fx). $$

Let us consider Barnes' type multiple $q$-Bernoulli polynomials.
For $w_1, w_2, \cdots, w_r \in \Bbb Z_p,$ and $ \delta_1,
\delta_2, \cdots,  \delta_r \in \Bbb Z,$ we define Barnes' type
multiple $q$-Bernoulli polynomials as follows:
$$ \aligned   & \beta_{n, q}^{(r)}(x \mid  w_1, \cdots, w_r : \delta_1, \cdots,  \delta_r ) \\
& = \underset{r \mbox{ times}}{\underbrace{ \int_{  \mathbb{Z}_p
}\cdots \int_{  \mathbb{Z}_p }
 }}    [w_1 x_1+ \cdots + w_r x_r +x ]_q^n   q^{ \sum_{j=1}^r(\delta_j-1)x_j}  d\mu_{q}(x_1)  \cdots d\mu_{q}(x_r).
 \endaligned  \eqno(30)
 $$
From (30), we can easily derive the following equation:
$$ \aligned   & \beta_{n, q}^{(r)}(x \mid  w_1, \cdots, w_r : \delta_1, \cdots,  \delta_r ) \\
& = \dfrac{1}{(1-q)^n} \sum_{l=0}^n \binom nl (-1)^l q^{lx}
\dfrac{(l w_1+ \delta_1)(l w_2+ \delta_2) \cdots (l
w_r+\delta_r)}{[l w_1+ \delta_1]_q [l w_2+ \delta_2]_q \cdots [l
w_r+\delta_r]_q }.
 \endaligned
 $$
Let $ \delta_r = \delta_1+r-1$.  Then we have
$$  \beta_{n, q}^{(r)}(x \mid  \underset{r \mbox{ times}}{\underbrace{ w_1 \cdots
w_1
 }}  : \delta_1,  \delta_1+1 \cdots,  \delta_1 +r-1)  = \dfrac{1}{(1-q)^n} \sum_{l=0}^n \binom nl (-1)^l q^{lx}
\dfrac{\binom{lw_1+\delta_1+r-1}{r}}
{\binom{lw_1+\delta_1+r-1}{r}_q } \dfrac{r!}{[r]_q!}.
  $$
Therefore, we obtain the following theorem.

\bigskip
{ \bf Theorem 11.} For $  w_1 \in \mathbb{Z}_p, r  \in \mathbb{N}$
and $\delta_1  \in \Bbb Z,$  we have
$$ \aligned  & \beta_{n, q}^{(r)}(x \mid  \underset{r \mbox{ times}}{\underbrace{ w_1 \cdots
w_1
 }}   : \delta_1,  \delta_1+1 \cdots,  \delta_1 +r-1) \\
& \quad  = \dfrac{1}{(1-q)^n} \sum_{l=0}^n \binom nl (-1)^l q^{lx}
\dfrac{\binom{lw_1+\delta_1+r-1}{r}}
{\binom{lw_1+\delta_1+r-1}{r}_q } \dfrac{r!}{[r]_q!}.
\endaligned
  $$

\end{document}